\documentclass{nyj}

\title{The canonical height of an algebraic 
point on an elliptic curve} 

\author{G. Everest}
\address{School of Mathematics, University of East Anglia,
Norwich NR4 7TJ, UK.} 
\email{g.everest@uea.ac.uk} 
\author{T. Ward}
\address{School of Mathematics, University of East Anglia,
Norwich NR4 7TJ, UK.} 
\email{t.ward@uea.ac.uk} 



\newtheorem{theorem}{Theorem}

\theoremstyle{definition}
\newtheorem{definition}[theorem]{Definition}
\newtheorem{example}[theorem]{Example}

\newtheorem{remark}[theorem]{Remark}

\begin{document} 

\begin{abstract}  
We use elliptic
divisibility sequences to describe
a method for estimating the global canonical height of an
algebraic point on an elliptic curve. This method requires 
almost no knowledge of the number field or the
curve, is simple to implement, and requires no
factorization.
The method is ideally suited to searching for algebraic points
with small height, in connection with the elliptic Lehmer problem. 
The accuracy of the method is discussed.
\end{abstract}
\maketitle
\tableofcontents


\section{Introduction}

Let $K$ denote an algebraic number field,
with ring of algebraic integers
$O_K$,
and $E$ an elliptic curve defined over $K$,
given by a generalized Weierstrass equation 
\begin{equation}
\label{geneq}
y^2+a_1xy+a_3y=x^3+a_2x^2+a_4x+a_6,
\end{equation}
with coefficients $a_1,\dots,a_6\in{O_K}$.
Let $Q=(x,y)$ denote a $K$-rational point of $E$, $Q\in E(K)$.
The {\sl global canonical height}
is a function $\hat h:E(K)\rightarrow {\mathbb R}$
with the properties:
\begin{enumerate}
\item $\hat h(Q)=0$ if and only if $Q$ is a torsion point of $E(K)$.
\item$\hat h(P+Q)+\hat h(P-Q)=2\hat h(P)+2\hat h(Q)$ for all $P,Q\in
E(K)$.
\end{enumerate}
The second property is known as the {\sl parallelogram law}.
The global canonical height is of fundamental importance in the
arithmetic
of elliptic curves, due in part to its functoriality.
The height appears in basic conjectures
such as Birch-Swinnerton-Dyer and there is a deep conjecture
known as the {\sl elliptic Lehmer problem}, see
\cite{MR92e:11062},
concerning lower bounds for the height.
Besides theoretical considerations however, it sometimes happens
that one really wishes to compute the value of the height
(for example, to compute the determinant of the height-regulator
matrix in searching for curves of large rank).

Silverman \cite{MR89d:11049} described an
algorithm for computing the global height, which can be made
arbitrarily accurate.
In the rational case, this is
implemented in Pari-GP (see \cite{pari-gp}). 
The  algorithm in \cite{MR89d:11049} requires the
discriminant of the curve to be completely factorized;
computing the height when the
discriminant cannot be factorized in reasonable time is
considered in \cite{MR97f:11040}. Silverman's method
requires much less than the full factorization of the
discriminant but still requires some factorization.
In
principle,
the method extends to the general algebraic
case, though there is currently no implementation of Silverman's algorithm in
the general case. When it is implemented, it is likely
to enjoy the same high accuracy it does in the rational case.
However there seems to be a small class of curves for which it is
vulnerable (see Examples \ref{silnodo}
and \ref{stillnodo} in Section \ref{examples}). 

Tate's definition of the global height gives a factorization-free
approach to computing the global height.
Let $M_K$ denote the set of valuations of $K$, each one
corresponding to an absolute value $\vert\cdot\vert_v$ (see
\cite{MR55:302} for background). 
For each valuation $v\in M_K$, let $K_v$ denote the
corresponding completion of $K$.
The {\sl naive height} $h(\alpha)$
of $\alpha\in K$ is
\begin{equation}\label{naiveheight}
h(\alpha)=\frac{1}{d}
\sum_{v\in M_K}\log \max \{1,\vert \alpha \vert_v\}.
\end{equation}
For a finite point $Q\in E(K)$,
set
$h(Q)=h(x(Q))$, and for $Q$ the point at infinity
set $h(Q)=0$.
Tate's definition of the global canonical height is
\begin{equation}\label{tatesdefinition}
\hat h(Q)={\frac12}\lim_{n\rightarrow \infty}4^{-n}h(2^nQ).
\end{equation}
Knowledge of the naive height is essentially
equivalent to knowledge of the minimal polynomial.

Tate's definition is not usually considered to be a very
useful method for actually computing the height.
In principle it is accurate: However,
it requires the computation of large integers and this 
not only slows
it down but makes high accuracy impossible in practice. 
On the other hand, 
it does always give an answer
because no factorization is needed.

The aim of this note is to exhibit an alternative 
factorization-free method for computing
the global height of an algebraic point on an elliptic curve
which stands somewhere between the two algorithms above.
Like the method in (\ref{tatesdefinition}), ours is extremely
simple, requiring almost no knowledge of the number field or
the elliptic curve and it does not require the curve to be in minimal
form.
However, our method gives more information than (\ref{tatesdefinition}) since
it also yields the archimedean and non-archimedean parts of the height
separately. (If the factorization of the discriminant is known
then it will give a complete decomposition of the global height
as a sum of local heights.) Our method can also be made much quicker.
It gives less accuracy 
than 
Silverman's algorithm but high accuracy is not always required.
In certain cases, our method can be used in tandem with
Silverman's algorithm: see Example \ref{tandem}
in Section \ref{examples}.

An example of a calculation which does not require great accuracy
is the search for algebraic points with small height. This
requires 
an accuracy of only 3 or 4
significant figures together with an easy way of handling algebraic number fields. 
Calculations such as these would shed light on the elliptic
Lehmer problem. In the classical Lehmer problem and its derivatives
(see \cite{MR2000e:11087}) there are many numerical examples.
Up to now, there is very little data for the elliptic Lehmer problem
outside the rational case. To illustrate our method, we give a
couple of examples of small height points found with an easy search: 
See Examples \ref{elp-3} and \ref{elp5} in
Section \ref{examples}.
Our method uses {\sl elliptic divisibility
sequences}, which are sequences associated to the division points
on the curves.
At the conclusion of the paper, we will make the point that our
methodology not only gives a simple way of handling elliptic curves
over algebraic number fields; it also throws up the possibility that
small height points might be found more efficiently by searching
for growth rates of elliptic divisibility sequences.

\section{Elliptic divisibility sequences}

The essential ingredient in the approach taken here is
the sequence of {\sl elliptic division polynomials}. For background on
see \cite{MR87g:11070} and \cite{MR96b:11074}.
\begin{definition}\label{beehive} With the notation of (\ref{geneq}),
define
\begin{eqnarray*}
b_2&=&a_1^2+4a_2,\\
b_4&=&2a_4+a_1a_3,\\
b_6&=&a_3^2+4a_6,\\
b_8&=&a_1^2a_6+4a_2a_6-a_1a_3a_4+a_2a_3^2-a_4^2.
\end{eqnarray*}
Define a sequence $(\psi_n)$ of polynomials in $O_K[x,y]$
as follows: $\psi_0=0, \psi_1=1$,
\begin{eqnarray*}
\psi_2&=&2y+a_1x+a_3,\\
\psi_3&=&3x^4+b_2x^3+3b_4x^2+3b_6x+b_8,\mbox{ and }\\
\psi_4&=&\psi_2(2x^6+b_2x^5+5b_4x^4+10b_6x^3+10b_8x^2+
(b_2b_8-b_4b_6)x+b_4b_8-b_6^2).
\end{eqnarray*}
Now define inductively for $n\ge 2$
\begin{eqnarray*}\label{oddandeven}
\psi_{2n+1}&=&\psi_{n+2}\psi_n^3-\psi_{n-1}\psi_{n+1}^3\mbox{ and }\\
\psi_{2n}\psi_2&=&\psi_n(\psi_{n+2}\psi_{n-1}^2-\psi_{n-2}\psi_{n+1}^2).
\end{eqnarray*}
\end{definition}
It is straightforward
to check that each $\psi_n\in O_K[x,y]$. It is known that
$\psi_n^2$ is a polynomial in $x$ alone having degree $n^2-1$
and leading coefficient $n^2$. The zeros of $\psi_n$
are the $x$-coordinates
of the points on $E$ with order dividing $n$.
Write $\psi_n(Q)$
for $\psi_n$ evaluated at the point $Q=(x,y)$. The sequence
$\psi_n(Q)$ is known as an {\sl elliptic divisibility sequence}:
Writing $u_n=\psi_n(Q)$ gives the elliptic recurrence relation
\begin{equation}\label{elldivseqrec}
u_{m+n}u_{m-n}=u_{m+1}u_{m-1}u_{n}^2-u_{n+1}u_{n-1}u_{m}^2
\end{equation}
for all $m \ge n\ge 0$.
These elliptic divisibility sequences were studied in an abstract
setting by Morgan Ward in a series of papers -- see \cite{MR9:332j}
for the details. Shipsey's thesis \cite{shipsey-thesis} contains
more recent applications of these sequences, which
satisfy the same recursion formul\ae{\ }as the division polynomials.
If $Q$ is
not a torsion point then the terms of the sequence $(\psi_n(Q))$ are
always non-zero.
The single relation (\ref{elldivseqrec})
gives
rise to the two relations 
\begin{equation}\label{oddrelation}
u_{2n+1}=u_{n+2}u_n^3-u_{n-1}u_{n+1}^3,\quad\mbox{ and }
\end{equation}
\begin{equation}\label{evenrelation}
u_{2n}u_2=u_{n+2}u_nu_{n-1}^2-u_{n}u_{n-2}u_{n+1}^2.
\end{equation}

For computational purposes, it is useful
to notice that
the two relations (\ref{oddrelation}) and (\ref{evenrelation}) 
can be subsumed 
into the single relation
$$
u_nu^{\vphantom2}_{\lfloor n/\lfloor
(n+1)/2\rfloor\rfloor}=
u^{\vphantom2}_{\lfloor(n+4)/2\rfloor}u^{\vphantom2}_{\lfloor n/2\rfloor}
u_{\lfloor(n-1)/2\rfloor}^2-
u^{\vphantom2}_{\lfloor(n+1)/2\rfloor}
u^{\vphantom2}_{\lfloor(n-3)/2\rfloor}u_{\lfloor(n+2)/2\rfloor}^2,
$$
where $\lfloor\cdot\rfloor$
denotes, as usual, the integer part.

Write
\begin{equation}\label{disc}
\Delta=-b_2^2b_8-8b_4^3-27b_6^2+9b_2b_4b_6 \in O_K
\end{equation}
for the discriminant of the curve $E$.
The valuations $v$ with $|\Delta|_v<1$
are precisely the valuations
corresponding to primes at which $E$ reduces to a singular curve.
Let $D=N_{K|{\mathbb Q}}(\Delta)$ and write $T$ for the set of rational primes
which divide $D$.
Given an algebraic integral point $Q\in E(K)$, let
$$E_n=|N_{K|{\mathbb Q}}(\psi_n(Q))| \mbox{ and } 
F_n=|E_n|\prod_{p\in T}|E_n|_p.
$$ 
Our method comes from the following theorem.

\begin{theorem}\label{mainformula}Let $Q$ denote an algebraic
integral point on $E(K)$. Then
\begin{equation}\label{maintheoremlimit}
\hat h(Q)=\frac1{d}\lim_{n\rightarrow \infty}\frac1{n^2}\log F_n.
\end{equation}
The total archimedean contribution is the limit
\begin{equation}\label{archcontribution}
h_{\infty}(Q)=\frac1d\lim_{n\rightarrow \infty}\frac1{n^2}\log E_n.
\end{equation}
\end{theorem}
The formula (\ref{maintheoremlimit}) is independent of the equation 
defining the curve.
It might appear that a factorization of the discriminant is
required but that is not so. Later
we discuss the practicalities of implementing the method.
The method extends to rational points provided one knows
the valuations at which the $x$-coordinate is not integral.
The
denominator can be cleared
to obtain an integral point on a curve
isomorphic to the starting curve, so the height is unchanged.
The proof of Theorem \ref{mainformula} follows in the next section.
It uses some detailed knowledge of local heights. For readers
interested
only in the application of the formula, the next section can be skipped.

\section{Local and global heights}

The global height is known to be expressible as a sum of
local heights, one for each element of $M_K$.  
There is a
function, continuous away from infinity, 
$\lambda_v:E({{\mathbb Q}}_v)\to{{\mathbb R}}$
which satisfies the {\sl local parallelogram law}
\begin{equation}\label{paralawlocal}
\lambda_v(P+Q)+\lambda_v(P-Q)=2\lambda_v(Q)+2\lambda_v(P)-\log
\vert x(Q)-x(P)\vert_v.
\end{equation}
Let $n_v=[K_v:{\mathbb Q}_w]/[K:{\mathbb Q}]$ denote the usual
local normalizing constants, where $v$
lies above $w$ on ${\mathbb Q}$. 
Then
\begin{equation}
\hat h(Q) = \sum_{v\in M_K}n_v\lambda_v(Q).
\end{equation}
The fundamental observation behind
our method is the {\sl elliptic Jensen formula}
from \cite{MR97e:11064}. If $G$ is a
compact group containing $Q$, with normalized Haar
measure $\mu_G$, then
\begin{equation}\label{elljenformula}
\lambda_v(Q)=2\int_G\log |x(P)-x(Q)|_vd\mu_G(P)
\end{equation}
by integrating and canceling three terms in (\ref{paralawlocal}).

If it is required that the expression $\lambda_p(Q)-
\frac{1}{2}\log\vert x(Q)\vert_p$ be bounded
as $Q\longrightarrow 0$, then there is only one
such map, the {\sl canonical local height}.
It is important to note that in \cite{MR96b:11074},
local heights are normalized to make them invariant
under isomorphisms. This involves adding a constant which
depends on the discriminant of $E$. The local heights
in \cite{MR96b:11074} satisfy a different form
of (\ref{paralawlocal}).

There are explicit formul{\ae}{\ }for each of the
local heights (see \cite{MR87g:11070} and
\cite{MR96b:11074}, or \cite{MR2000g:11050} for
an alternative approach).
For non-archimedean valuations $v$ where $Q$ has good reduction,
\begin{equation}
\lambda_v(Q)=\textstyle\frac{1}{2}\log \max \{1,\vert x(Q)\vert_v\}.
\end{equation}
Notice in particular that if $x(Q)$ is integral at $v$
and $Q$ has good reduction at $v$ then $\lambda_v(Q)=0$.
The bad reduction case is more involved but we need to
deal only with {\sl split multiplicative
reduction} (see \cite[p. 362]{MR96b:11074} for details on this).
This is because we may pass to an extension field where the
reduction becomes of this type --- the local height is
functorial in the sense that it respects this passage. 
In the split multiplicative case, the points on the
curve are isomorphic to the points on the
Tate curve ${K}_v^{*}/q^{{\mathbb Z}}$,
where $q\in K_v^{*}$ has $\vert q\vert_v<1$. The
explicit formul{\ae}{\ }for the $x$ and $y$ coordinates of a
non-identity point are given in terms of the uniformizing parameter
$u\in K_v^{*}$ by
\begin{eqnarray*}
x&=&x_u=\sum_{n\in{{\mathbb Z}}}\frac{q^nu}{(1-q^nu)^2}-
2\sum_{n\ge1}\frac{nq^n}{(1-q^n)^2},\\
y&=&y_u=\sum_{n\in{{\mathbb Z}}}\frac{q^{2n}u^2}{(1-q^nu)^3}+
\sum_{n\ge1}\frac{nq^n}{(1-q^n)^2}.
\end{eqnarray*}
It is clear that $x_u=x_{uq}$ and $x_u=x_{u^{-1}}$.
Suppose $Q$ corresponds to the point $u\in K_v^{*}$
and assume, by invariance under multiplication by $q$,
that $u$ lies in the fundamental domain
$\{u\mid\vert q\vert_v<\vert u\vert_v
\le 1\}$.
Then (by \cite{MR2000g:11050} or
\cite{MR96b:11074}), writing $\rho=\log \vert u\vert_v/\log \vert q\vert_v$,
$$
\lambda_v(Q)=\left\{
\begin{array}{ll}
-\log\vert 1-u\vert_v&\mbox{if }\vert u\vert_v=1,\\
\frac12(\rho-\rho^2)
\log\vert q \vert_v&\mbox{if }\vert u\vert_v<1.
\end{array}
\right.
$$
Notice that for $\vert u\vert_v=1$,
the local height is non-negative, while if $\vert u\vert_v<1$
the local height is negative.

\begin{theorem}\label{heightisthelimit} Let $Q$ denote a non-torsion integral
point. Suppose $v|\infty$ or $v$ corresponds to a prime of singular
reduction.
In the latter case, assume equation  {\rm(\ref{geneq})} is in minimal
form.
Then there are positive constants $A$ and $B<2$ such that 
\begin{equation}\label{localheightisthelimit}
\frac{1}{n^2}\log |\psi_n(Q)|_v = \lambda_v(Q)+
\left\{
\begin{array}{cl}
O((\log n)^A/n^2)&\mbox{ if } v|\infty, \\ 
O(1/n^B)&\mbox{ otherwise. }\end{array}\right.
\end{equation}
\end{theorem}

\begin{proof}
If $v|\infty$,
we claim first that
\begin{equation}
\label{firststepininfinitecase}
\lim_{n\to\infty}n^{-2}\log\vert\psi_n(x(Q))\vert_v=\lambda_{\infty}(Q).
\end{equation}
Formula (\ref{firststepininfinitecase}) was proved in the rational
case
in
\cite[Theorem 6.18]{MR2000e:11087}; the proof is
sketched here in the general case. The height is functorial
in the sense that it respects field extensions. Thus we
may assume $v$ corresponds to an embedding of $K$ into ${\mathbb C}$.
Take $G=E({{\mathbb C}})$ in the elliptic Jensen
formula (\ref{elljenformula}).
The points of $n$-torsion are dense
and uniformly distributed in $E({{\mathbb C}})$ as $n\to\infty$,
so the limit sum over the torsion points will tend to the
integral when the integrand is continuous. Note that the torsion
points occur in pairs usually. Working with $\psi_n(Q)$ they
only occur with multiplicity 1, hence the formula differs
from the usual elliptic Jensen formula in this respect.
The only potential problem arises from torsion points
close to $Q$: by \cite{MR98f:11078}, for $x=x(Q)$
with $nQ=0$,
$\vert x-x(Q)\vert_v>n^{-C}$ for
some $C>0$ which depends on $E$ and $Q$
only. This inequality is enough to imply
that the Riemann sum given by the $n$-torsion
points for $\log\vert x-x(Q)\vert_v$ converges,
which gives (\ref{firststepininfinitecase}), and
the explicit error term gives the estimate
in (\ref{localheightisthelimit}).

Assume now that $v$ is non-archimedean, corresponding
to a prime of singular reduction. Let $\Omega_v$ denote
any complete, algebraically closed field containing $K_v$.
Assume $Q$ is integral, $\vert x(Q)\vert_v\le 1$.
Now use the parametrisation of the curve
described before.
The points of order dividing $n$ on the Tate curve
are precisely those of the form $\zeta^iq^{j/n}$,
$1\le i,j\le n$, where $\zeta\in\Omega_v$ denotes a fixed,
primitive $n$th root of unity in $\Omega_v$.
We claim that
\begin{equation}\label{andthereisnoplace}
\lim_{n\to\infty}n^{-2}\log\vert\psi_n(x(Q))\vert_v=\lambda_{v}(Q).
\end{equation}
Let $G$ denote the closure of the torsion points:
$G$ is not compact, so the $v$-adic elliptic Jensen
formula cannot be used. Instead we use a variant of
the Shnirelman integral: for $f:E(\Omega_v)\to{{\mathbb R}}$ define
the elliptic Shnirelman integral to be
$$
\int_Gf(Q)\mbox{d}Q=
\lim_{n\to\infty}n^{-2}
\sum_{n\tau=0}f(\tau)
$$
whenever the limit exists.

We claim firstly that for any $P\in E({{\mathbb Q}}_p)$, the Shnirelman
integral
\begin{equation}\label{imgoingto}
\int_G\lambda_v(P+Q)\mbox{d}Q
=S(E)\mbox{ exists and is independent of $P$.}
\end{equation}
First assume that $P$ is the identity.
Using the explicit formula for the local height
gives
\begin{equation}\label{playasongforme}
-n^{-2}\sum_{i=1}^{n-1}\log\vert1-\zeta^i\vert_v
-n^{-2}\sum_{i=0}^{n-1}\sum_{j=1}^{n-1}
\frac{k}{2}\left(\frac{j}{n}-\left(\frac{j}{n}\right)^2\right).
\end{equation}
The first sum is bounded by $\log\vert n\vert_v/n$, which
vanishes in the limit; the second sum converges to
$-\frac{k}{12}$. For the general case, let $P$ correspond to the point
$u$ on the multiplicative Tate curve.
If for some large $n$ no $j$ has $\vert q^{j/n}u\vert_v=1$
then the analogous sum to (\ref{playasongforme})
is close to $-\frac{k}{12}$ by the same argument.
Assume therefore that there is a $j$ with this property.
Then the first sum in (\ref{playasongforme}) is replaced
by
\begin{equation}\label{imnotsleepy}
-n^{-2}\sum_{i=0}^{n-1}\log\vert
1-q^{j/n}u\zeta^i\vert_v
-n^{-2}\log\vert1-(q^{r}u)^n\vert_v,
\end{equation}
where $r=j/n$ only depends on $u$.
By $v$-adic elliptic transcendence theory (see
\cite{MR98f:11078}), there is
a lower bound for $\log\vert1-(q^ru)^n\vert_v$ of the
form $-(\log n)^A$, where $A$ depends on $E$ and $u=u(P)$
only. It follows that the first sum vanishes in
the limit as before. The second sum in
(\ref{playasongforme}) is simply rearranged under
rotation by $u$, so converges to $-\frac{k}{12}$ as before.
This proves (\ref{imgoingto}).

The claimed limit (\ref{andthereisnoplace})
now follows by taking the elliptic Shnirelman
integral of both sides of the parallelogram law
(\ref{paralawlocal}) and noting that we count torsion
points in pairs.  Equation (\ref{imgoingto})
shows that three terms cancel to leave the required limit.
The error term in (\ref{localheightisthelimit}) comes from
the lower bound used above.
\end{proof}

These estimates are enough to prove the main formula.

\begin{proof} (of Theorem \ref{mainformula})
It will be convenient to use normalized heights,
so define
$$
\nu_v(Q)=\lambda_v(Q)-{\textstyle\frac{1}{12}}\log |\Delta|_v.
$$
Then $\nu_v$ is invariant under
isomorphism (see \cite{MR96b:11074}). By the product formula,
$$\hat h(Q)=\sum_v n_v\nu_v(Q)=\sum_v n_v\lambda_v(Q).
$$
Also, by Theorem \ref{heightisthelimit},
\begin{equation} \lim_{n\rightarrow \infty}\frac1{n^2} 
\log |\psi_n(Q)\Delta^{{-n^2}/{12}}|_v=\nu_v(Q).
\end{equation}  
For any $\alpha \in K$, $|N_{K|{\mathbb Q}}(\alpha)|
=\prod_{v|\infty}|\alpha|_v.$ Therefore, using the product
formula again,
$$\log |F_n|=\sum_{v|\infty}\log |\psi_n(Q)\Delta^{{-n^2}/{12}}|_v
+\sum_{|\Delta|_v<1}|\psi_n(Q)\Delta^{{-n^2}/{12}}|_v.
$$
The reason for introducing the factor
$\Delta^{{-n^2}/{12}}$
is to take account
of the possibility  that the equation 
(\ref{geneq}) is
not in minimal form at some  non-archimedean $v$ 
corresponding to a prime of singular
reduction. The change of coordinates to put the
equation into minimal form is an isomorphism, so it leaves
the local height $\nu_v(Q)$ invariant. Now
Theorem \ref{mainformula} follows directly from Theorem
\ref{heightisthelimit}.  
\end{proof}

\section{Examples}
\label{examples}

It appears as though we need to factorize
$D=N_{K|{\mathbb Q}}(\Delta)$
in order to apply Theorem \ref{mainformula}. However,
Theorem \ref{heightisthelimit} says that for a
prime $p\in T$, $|E_n|_p$ is approximately $l^{n^2}$
where $l$ is the total contribution to the height
from the valuations which extend $\vert\cdot\vert_p$.
Therefore, asymptotically, it suffices to compute
the gcd of $E_n$ with a suitably high power of $D$.
Since the local height is $t\log \vert \Delta \vert_v$
for some $0\le t\le 1$, the power of $D$ can be $n^2$.
This is likely to be a huge number and there are ways to
avoid making this computation. In practice,
it is often
sufficient to find the gcd of $E_n$
and $E_{n+1}$. In other words:
\begin{equation}\label{bogstandard}
\hat h(Q)=\frac1{d}\lim_{n\rightarrow \infty} 
\frac1{n^2}\log \left(\frac{E_n}{\gcd(E_n,E_{n+1})}\right).
\end{equation}
In the last section of the paper, we will discuss other ways to
speed up the calculations.

The following examples were
calculated using Pari-GP, see \cite{pari-gp},
simply applying the basic formula (\ref{bogstandard}). 
In the main we have only exhibited 
calculations which
were executed within a few seconds at most.
We begin
by applying our method to examples in the literature --
the first two examples come from \cite{MR89d:11049}. 
\begin{example}\label{silver1}Let the curve be
$$E:y^2+y=x^3-x^2,$$
the field $K={\mathbb Q}(\sqrt{-2})$, and $Q=(2+\sqrt{-2},1+2\sqrt{-2})$.
Taking $n=100$ gives $\hat h(Q)\sim .45744\dots$ to be compared with
Silverman's accurate value of $.45754\dots$. When $n=200$, we obtain the better
approximation $\hat h(Q)\sim .45753\dots$.
\end{example}
\begin{example}\label{silver2}Let $K={\mathbb Q}(i)$, 
let the curve be
$$E:y^2+4y=x^3+6ix,$$
and $Q=(0,0)$. Taking 
$n=200$ gives  $\hat h(Q) \sim .33688\dots$ 
to be compared with
Silverman's accurate value of $.33689\dots$ The archimedean
height is $\sim .51016\dots$
\end{example}

The next example illustrates that the curve does not need to
be in minimal form for the method to work.
\begin{example}\label{notinmin}Let the curve be
$$E:y^2=x^3-16x+16,$$
and let $Q=(0,4)$. Taking $n=150$ gives a value $\hat h(Q)
\sim .02549\dots$ with a value $\sim .7186\dots$ for the archimedean
component. The calculation speeds up if we notice that $E$ is
isomorphic to the curve
$y^2+y=x^3-x,
$
with $Q$ mapping to $P=(0,0)$ under the isomorphism.
Taking $n=150$ gives 
$\hat h(P)=\hat h(Q) \sim .02555\dots$ which is more
accurate, and quicker, due to the slower growth rate
of the sequence $E_n$.
\end{example}

The next examples are manufactured to highlight one of the
strengths of our approach: It always gives an answer
even if a tricky factorization appears to be necessary.
Silverman's approach in \cite{MR97f:11040} computes all
the local heights then sums these to give the global
height. To compute a local non-archimedean height, the
curve needs to be in minimal form for that valuation.
If the factorization of $\Delta$ is known then
the curve can easily be rendered in minimal form for each
valuation corresponding to the prime factors of $\Delta$.
Even if the factorization is not known, it is usually
possible to proceed.
With our earlier notation,
define
$$c_4=b_2^2-24b_4 \mbox{ and } c_6=-b_2^3+36b_2b_4-216b_6.
$$
In \cite{MR97f:11040}, working over $\mathbb Q$,
Silverman shows that if the
factorization of $c=$
gcd$(c_4,c_6)$ is known then the curve can be put in global
minimal
form so the local heights can all be computed. Over a number
field with class number greater than 1, a global minimal
equation will not always exist. Presumably the same
kind of argument would work nonetheless. 
Therefore, the next example is chosen
to highlight a potential difficulty:
$c$ may have a large gcd with the discriminant.
In this case, factorizing
$c$ is not much easier than factorizing
the discriminant.

\begin{example}\label{silnodo} Let $K={\mathbb Q}$ and let $m\in {\mathbb N}$
denote an integer that is not factorizable in reasonable
time.
Consider the curve
$$E:y^2=x^3+mx+m^2.$$
Let $Q$ denote the point $(0,m)\in E(K)$.
For this curve, $m|c$ and Silverman's algorithm 
now requires auxiliary arguments (see Remark \ref{finalremark}
below).
Let $m=pq$
where $p$ and $q$ denote the next primes after $10^{30}$
and $10^{40}$. With $n={50}$, within a minute, our method gave 
$\hat h(Q)\sim13.657\dots$
with an archimedean height $\sim 53.936\dots$.
We also
used a floating point for the archimedean contribution: with $n=300$
we obtained $\sim 53.956\dots$.
The Pari-GP routine for computing heights
returned a warning that the calculation would take several hours.
This is all due to the difficulty of factorizing $m$.
\end{example}

\begin{remark}\label{finalremark}
The referee pointed out to us that Silverman's method can be made
to work in this example because 
it can be checked that no 4th power of a prime divides $m$.
\end{remark}

The next example shows how our method can be used in tandem
with Silverman's
algorithm.
\begin{example}\label{tandem} With $E$ as in the previous example, 
let $Q$ denote an algebraic point with $x(Q)=1$. 
Even a small value of $n$ shows the total
non-archimedean contribution is zero. Thus one may revert
immediately to a general algebraic version of
Silverman's method to obtain a very accurate
value for the global height, which is entirely concentrated
at the archimedean valuation. Using our method, with $n=300$ and 
with floating point arithmetic
on the two archimedean valuations, we obtained the
value $\sim 53.956\dots$ for the total archimedean contribution. 
Note the value is close to the previous example -- this is
no real surprise, as the archimedean heights are continuous.
\end{example}

Our next example is an algebraic version of Example \ref{silnodo}.

\begin{example}\label{stillnodo} Let $f(x)=x^{17}
+x+996$ and let 
$K={\mathbb Q}(\rho)$ where $\rho$ denotes any root
of $f(x)$. Let $\theta=1-1728\rho^2$, and consider the curve 
$$E:y^2=x^3+\theta x+\theta^2.
$$
Let $Q$ denote the point $(0,\theta)\in E(K)$. With $n=35$,
in under one minute our method gives $\hat h(Q)\sim 15.595\dots$. 
The archimedean height is $\sim 50.732\dots$. As in the previous
example,
$\theta|c$.
It took Pari-gp 
30 minutes to find the factorization
$$C=
11978293086538309\times 904414027740749856394559037844972335934195571
$$
of $C=|N_{K|{\mathbb Q}}(\theta)|$;
it would have taken at least as long to factorize 
the ideal $(c)$.
\end{example}

Finally, we give two examples of small height points over algebraic
number fields. Our method is simple to apply and can be used to
search for small height points in connection with the elliptic Lehmer
problem. There is very little data associated
with this problem beyond the rational
case. We hope our paper might inspire an attempt to
gather some data.

\begin{example}\label{elp-3}Let $w$ denote a non-trivial cube root
of unity and $K={\mathbb Q}(w)$. Let $E$ be the elliptic curve
$$y^2=x^3-243x+3726+10368w.
$$
The point $Q=(3-12w,-108w^2)$ has global height $\hat h(Q) \sim .01032\dots$.
This was found taking $n=512=2^9$ and using Shipsey's algorithm
from the next section.
Although the coefficients of the curve might seem large, this
example arises from a simple elliptic divisibility sequence. Starting
from the sequence $0,1,1+w,1+w,1+w,\dots$ we used Morgan Ward's
formul{\ae}{\ }(see \cite[p. 50]{MR9:332j}) to
obtain a point on a curve
with coefficients in $K$ whose denominators can be cleared to give $E$
as above. 
\end{example} 

\begin{example}\label{elp5}Let $u=(1+\sqrt 5)/2$ and $K={\mathbb Q}(u)$.
The curve $E$ is
$$y^2=x^3+(-2214+1215u)x+40878-23328u
$$
and the point is $Q=(3-9u,108-108u)$. Taking $n=512$ as before
gives $\hat h(Q) \sim .00971\dots$. This example came from the
elliptic divisibility sequence which begins in the modest way 
$0,1,1-u,-2+u,5-3u,\dots$ Inverting this sequence gives a point
on a curve over $K$ and clearing the denominators gives $E$
as above.
\end{example}

Two comments needs to be made about these examples. Firstly, although
these heights are small, no records have been broken. The elliptic
Lehmer problem predicts a lower bound for $d\hat h(Q)$ where $d$
is the degree of the number field. Multiplying both the above by 2
shows these values are not smaller than the height ($\sim .01028$)
of the rational
point $Q=(13,33)$ on the curve $y^2+xy+y=x^3-x^2-48x+147$, which
appears in \cite[p. 480]{MR96b:11074}. Secondly, these examples hint
at an interesting possibility concerning the search for small
height points. Perhaps restricting to elliptic divisibility
sequences represents an efficiency gain in the sense that
small height points will arise from sequences whose first
few terms are arithmetically simple.

\section{Accuracy}
In (\ref{localheightisthelimit}), the error term is estimated
using methods from elliptic transcendence theory. In \cite{eds},
we investigated the error in practice and found it to be about
$O(1/n^2)$, even for quite modest values of $n$. 
For small values of $n$, the values of $E_n$ can be computed
easily using Pari-GP.
Several options for achieving greater accuracy are listed below.
However, we stress again that there are certain physical limits to this
method which go beyond computational considerations:
Accuracy of 80 significant figures would involve computing
a number with approximately $10^{40}$ decimal digits.
Even storing such numbers is beyond the capabilities of
any computer.

1. The archimedean and non-archimedean contributions
can be computed separately and this allows the computations
to be speeded up. For the archimedean contribution, we can
use floating point arithmetic which greatly enhances the
speed. For the non-archimedean contribution, we only have
to keep a running total of the gcd so big integer arithmetic
can be avoided. If the factorization of the discriminant is
known then p-adic arithmetic may be used.

2. Since the computation of the height involves big numbers,
it is useful to use a package which allows these to be
handled efficiently.
We are grateful to John Cannon for implementing our algorithm
in Magma \cite{magma} which gave greater accuracy.

3. Memory is clearly an issue with the method we are describing
since it involves the calculation of huge numbers. Storage
can be maximized by computing $E_n$ for special $n$, without
needing to know all $E_m$ for $m<n$. Shipsey \cite{shipsey-thesis}
gives an algorithm that computes $E_{n}$ in $O(\log n)$
arithmetic operations. Note the distinction between arithmetic
operations and bit operations: By arithmetic operation is meant
one of the familiar operations of adding or multiplying.
The special case where $n=2^N$ is especially easy to
implement and we describe it below. We are grateful
to Rachel Shipsey for her permission to include it here.

Now follows Shipsey's algorithm for computing 
$E_n$ when $n=2^N$:
Given $Q$ and $E$, find $\psi_i(Q)$ for $i=2,3,\dots,7$
using the formulae given before. Let
$$
T_1=1,
U_1=\psi_2(Q),
V_1=\psi_3(Q),
W_1=\psi_4(Q),
X_1=\psi_5(Q),
Y_1=\psi_6(Q),
Z_1=\psi_7(Q),
$$
and then inductively
\begin{eqnarray*}
T_{n+1}&=&W_nU_n^3-V_n^3T_n,\\
U_{n+1}&=&(V_n/\psi_2(Q))(X_nU_n^2-T_nW_n^2),\\
V_{n+1}&=&X_nV_n^3-W_n^3U_n,\\
W_{n+1}&=&(W_n/\psi_2(Q))(Y_nV_n^2-U_nX_n^2),\\
X_{n+1}&=&Y_nW_n^3-X_n^3V_n,\\
Y_{n+1}&=&(X_n/\psi_2(Q))(Z_nW_n^2-V_nY_n^2),\\
Z_{n+1}&=&Z_nX_n^3-Y_n^3W_n.
\end{eqnarray*}
After $N-2$ iterations the value of $W$ is
$\psi_{n}(Q)$, and $E_{n}=|N_{K|{\mathbb Q}}(\psi_{n}(Q))|$.

Computing $E_n$ requires $O(\log n)$
arithmetic operations. The operations required for
(\ref{tatesdefinition}) satisfy the same bound. However,
our method can be speeded up in two ways. Firstly, by using floating 
point arithmetic
for the archimedean contribution. Secondly, the homogeneity
of the formul{\ae}{\ }make it possible to keep a running
total for the gcd computation, yielding the non-archimedean
contribution. By successively factoring out the gcd, the calculations
proceed with smaller integers, making the method much faster.


\begin{thebibliography}{War48}

\bibitem[Dav95]{MR98f:11078}
Sinnou David, {\em Minorations de formes lin\'eaires de logarithmes
  elliptiques}, M\'em. Soc. Math. France (N.S.) (1995), no.~62, iv+143.

\bibitem[EF96]{MR97e:11064}
G.~R. Everest and Br{\'\i}d~N{\'\i} Fhlath{\'u}in, {\em The elliptic {M}ahler
  measure}, Math. Proc. Cambridge Philos. Soc. {\bf 120} (1996), no.~1, 13--25.

\bibitem[Eve99]{MR2000g:11050}
Graham Everest, {\em Explicit local heights}, New York J. Math. {\bf 5} (1999),
  115--120 (electronic).

\bibitem[EEW]{eds}
M.~Einsiedler, G.~Everest and T.~Ward.
{\em Computational aspects of elliptic divisibility sequences},
Pre-print.

\bibitem[EW99]{MR2000e:11087}
Graham Everest and Thomas Ward, {\em Heights of Polynomials and Entropy in
  Algebraic Dynamics}, Springer-Verlag London Ltd., London, 1999.

\bibitem[HS90]{MR92e:11062}
Marc Hindry and Joseph~H. Silverman, {\em On {L}ehmer's conjecture for elliptic
  curves}, S\'eminaire de Th\'eorie des Nombres, Paris 1988--1989, Birkh\"auser
  Boston, Boston, MA, 1990, pp.~103--116.

\bibitem[GP]{pari-gp}
{\textsc{Pari-GP}}, {\tt http://www.parigp-home.de}.

\bibitem[Mag]{magma}
{\textsc{Magma}}, {\tt http://www.maths.usyd.edu.au:8000/u/magma}.

\bibitem[Shi00]{shipsey-thesis}
Rachel~Shipsey.
\newblock {\em Elliptic Divisibility Sequences}.
\newblock PhD thesis, University of London (Goldsmiths), 2000.

\bibitem[Sil86]{MR87g:11070}
Joseph~H. Silverman, {\em The Arithmetic of Elliptic Curves}, Springer-Verlag,
  New York, 1986.

\bibitem[Sil88]{MR89d:11049}
Joseph~H. Silverman, {\em Computing heights on elliptic curves}, Math. Comp.
  {\bf 51} (1988), no.~183, 339--358.

\bibitem[Sil94]{MR96b:11074}
Joseph~H. Silverman, {\em Advanced Topics in the Arithmetic of Elliptic
  Curves}, Springer-Verlag, New York, 1994.

\bibitem[Sil97]{MR97f:11040}
Joseph~H. Silverman, {\em Computing canonical heights with little (or no)
  factorization}, Math. Comp. {\bf 66} (1997), no.~218, 787--805.

\bibitem[War48]{MR9:332j}
Morgan Ward, {\em Memoir on elliptic divisibility sequences}, Amer. J. Math.
  {\bf 70} (1948), 31--74.

\bibitem[Wei74]{MR55:302}
Andr{\'e} Weil, {\em Basic Number Theory}, third ed., Springer-Verlag, New
  York, 1974, Die Grundlehren der Mathematischen Wissenschaften, Band 144.

\end{thebibliography}

\ifx\undefined\bysame
\newcommand{\bysame}{\leavevmode\hbox to3em{\hrulefill}\,}
\fi

\end{document}